\documentclass[amstex, 11pt]{amsart}
\def\depth{\mathop{\rm depth}}

\newcommand{\sta}{\stackrel}

\def\C{{\mathcal C}}
\def\F{{\mathcal F}}
\def\ia{{\mathfrak a}}

\def\m{{\mathfrak m}}
\def\Fi{{\mathcal F}_{I_1;I_2}}
\def\oFi{{\mathcal F}_{\overline{I_1};\overline{I_2}}}

\def\F2i{{\mathcal F}_{I_2}}

\def\Cnn{C({\bf x}_{k}, {\F2i}, (0, n))}

\def\C1n{  C({\bf x}_{k}, {\Fi}, (1, n))}

\def\Cd11{  C({\bf x}_{d}, {\Fi}, (1, 1))}

\def\Cdnn{  C({\bf x}_{d}, {\Fi}, (1, n))}
\def\Ddnn{  D({\bf x}_{d}, {\Fi}, (1, n))}

\def\Cdn0{  C({\bf x}_{d}, {\Fi}, (1, 0))}
\def\Cn1{C({\bf x}_{k}, {\Fi}, (1, n-1))}
\def\Cnk1{C({\bf x}_{k-1}, {\Fi}, (1, n))}

\def\Cn1k1{C({\bf x}_{k-1}, {\Fi}, (1, n-1))}

\def\D1n{D({\bf x}_{k}, {\Fi}, (1, n))}
\def\Dd11{D({\bf x}_{d}, {\Fi}, (1, 1))}
\def\Ddn0{D({\bf x}_{d}, {\Fi}, (1, 0))}
\def\Dn1{D({\bf x}_{k}, {\Fi}, (1, n-1))}
\def\Dnk1{D({\bf x}_{k-1}, {\Fi}, (1, n))}

\def\x{{\bf x}}

\newtheorem{thm}{Theorem}[section]
\newtheorem{lemma}[thm]{Lemma}
\newtheorem{cor}[thm]{Corollary}
\newtheorem{pro}[thm]{Proposition}
\newtheorem{example}[thm]{Example}
\newtheorem{remark}[thm]{Remark}
\newtheorem{defn}[thm]{Definition}
\newtheorem{notn}[thm]{Notation}
\newtheorem{blank}[thm]{}

\newcommand{\bt}{\begin{thm}}
\newcommand{\et}{\end{thm}}
\newcommand{\blem}{\begin{lemma}}
\newcommand{\elem}{\end{lemma}}
\newcommand{\bco}{\begin{cor}}
\newcommand{\eco}{\end{cor}}
\newcommand{\bp}{\begin{pro}}
\newcommand{\ep}{\end{pro}}
\newcommand{\bex}{\begin{example}}
\newcommand{\eex}{\end{example}}
\newcommand{\brm}{\begin{remark}}
\newcommand{\erm}{\end{remark}}
\newcommand{\bdefn}{\begin{defn}}
\newcommand{\edefn}{\end{defn}}
\newcommand{\bnot}{\begin{notn}}
\newcommand{\enot}{\end{notn}}
\newcommand{\bblank}{\begin{blank}}
\newcommand{\eblank}{\end{blank}}

\def\nno{\nonumber}

\newcommand{\bib}{\bibitem}
\newcommand{\f}{\frac}
\newcommand{\sms}{\setminus}
\newcommand{\beqn}{\begin{eqnarray*}}
\newcommand{\eeqn}{\end{eqnarray*}}
\newcommand{\beq}{\begin{eqnarray}}
\newcommand{\eeq}{\end{eqnarray}}
\newcommand{\been}{\begin{enumerate}}
\newcommand{\eeen}{\end{enumerate}}
\newcommand{\lrar}{\longrightarrow}
\newcommand{\rar}{\rightarrow}

\newcommand{\Z}{\mathbb Z}

\newcommand{\limn}{\underset{\underset{n}{\longrightarrow}}{\lim}}

\newcommand{\mptm}{\fontfamily{ptm}\fontsize{11}{13}\selectfont}
\newcommand{\mptmcomx}{\fontfamily{ptm}\fontsize{9}{11}\selectfont}

\newcommand{\ncent}{\fontfamily{pcr}\fontseries{b}\fontsize{12}{15}\selectfont}

\begin{document}

\title[Fiber cones]{
 On the Homology and   fiber cone of ideals}

\author{  Clare D'Cruz}
\address{Chennai Mathematical Institute, 
Plot H1, SIPCOT IT Park, 
Siruseri,
Paddur PO, 
Siruseri, Kelambakkam 603 103, 
India}
\email{clare@cmi.ac.in}
\mptm
\maketitle
\begin{abstract}
In this paper we give a unified approach for several results concerning  the fiber cone. Our novel ideal is to use the   complex 
$\C1n$. We improve earlier  results  obtained by several researchers and get  some new results. We  give a more general definition of  ideals of minimal multiplicity and  of ideals of almost minimal multiplicity. We also  compute the Hilbert series of the fiber cone for these ideals. 
\end{abstract}

\section{Introduction}
Throughout this paper we will assume that $(R,\m )$ is a local ring of positive dimension $d$ and  infinite residue field. The associated graded ring $G(I) := \oplus_{n \geq 0} I^n/ I^{n+1}$ has been investigated in detail by several researchers. In the last two decades  the fiber cone $F(I):= \oplus_{n \geq 0} I^n/ \m I^n$ has been of interest. 
  Let $I_1$ and $I_2$ be ideals in $(R, \m)$. We call  $F_{I_1}(I_2) := \oplus_{n   \geq 0} I_2^n/ I_1 I_2^n$ the fiber cone of $I_2$ with respect to $I_1$. 
  Since $G(I)= F_{I}(I)$,  it is of interest to know  how the properties of these two rings are related. Several recent papers on the fiber cone do imply that it is possible to extend the results known for the associated graded ring to the fiber one. 
  
  We begin by recalling a few results on  $G(\m)$. Let $(R,\m)$ be a Cohen-Macaulay ring. Sally  showed that if $\m$ is an ideal with   minimal multiplicity,  then $G(\m)$ is Cohen-Macaulay and  the corresponding Hilbert function $\ell(\m^n / \m^{n+1})$ can be explicitly described
  (\cite{sally1}).
  She conjectured that if $\m$ is an ideal with almost minimal multiplicity, then  $G(\m)$ has almost maximal depth
(\cite{sally2}).  Her conjecture was settled independently by Rossi and Valla in \cite{rossi-valla}   and by Wang in \cite{wang}.  Sally's work has been generalized in various directions.

Goto gave a more general definition of  ideals of minimal multiplicity \cite{goto}. 
Inspired by his work, Jayanthan and Verma defined ideals of minimal multiplicity and ideals of almost minimal multiplicity in the case when $I_1$ and $I_2$ are $\m$-primary ideals satisfying 
$I_2 \subseteq I_1$ 
(\cite{jay-verma}, \cite{jay-verma-almm}). They studied the fiber cone of these ideals in great detail. They  generalized Sally's conjecture to ideals of almost minimal multiplicity and showed that if  the depth of $G(I_2)$ is at least  $d-2$, then the  depth of 
$F_{I_1}(I_2)$ is atleast  $d-1$. 
In this paper, we define ideals of minimal multiplicity and ideals of almost minimal multiplicity for any two $\m$-primary ideals (Definition~\ref{defn-mm-amm}, Lemma~\ref{min-mult-defn}).

We  take a new approach in this paper.
  Let  $1 \leq k \leq d$ and let $\x_{k}:= x_1, \ldots, x_k$  a system of parameters in $I_2$. 
 The complex $\Cnn$
has been studied in    \cite{anna} and \cite{tom-huc} in connection with the properties of the associated graded ring $G(I)$.  
 They showed that   vanishing of the complex $C(\x_{k}, {\Fi}, (i, n))$  determines 
 $\depth ((\x_k)^{\star}, G(I))$, where $x_i^{\star}$ denotes the image of $x_i$ in $I/I^2$. 
    Associated to the $I_2$-filtration $\Fi = \{I_1I_2^n\}_{n \in \Z}$, 
  we have the 
  complex $\C1n$ (defined in Section~\ref{complex}). The corresponding graded ring 
  $G_{I_1}(I_2)
  :=  R/I_1 \oplus \oplus_{n \geq 1} ({\Fi})_{n-1}/ ({\Fi})_{n}$  has been studied in \cite{rossi-valla-2}.
 Note that $G(I) = G_{R}(I)[1]$. 
 When $(\x_k) \subseteq I_1$, we consider the truncated complex 
$D(\x_{k}, {\Fi}, (i, n))$ 
and one can verify  that $D(\x_{k}, {\Fi}, (0, n))=C(\x_{k}, {\Fi}, (0, n))$.
 For any element $x \in I_2$, let  ${x^{o}}$ denote the image of $x$ in $I_2/ I_1 I_2$. 
In this paper, we use the complexes $D(\x_{k}, {\Fi}, (i, n))$ (i=0,1) and the  Koszul complex
  $K({\x_{k}}^o, F_{I_1}(I_2))(n)$ to investigate the relation between the properties of the three graded rings 
 $G(I_2)$, $G_{I_1}(I_2)$ and $F_{I_1}(I_2)$. As a consequence,  
we  obtain interesting information on the fiber cone.

Huneke's fundamental lemma (\cite[Lemma~2.4]{huneke})  was extended to the filtration $\Fi$ for $\m$-primary ideals  $I_2 \subseteq I_1$ in a two dimensional Cohen-Macaulay local ring 
 (\cite[Proposition~2.5]{jay-verma-almm}). 
 In this  paper, using the complex
 $D(\x_{d}, {\Fi}, (1, n))$ we    extend  this  result to   any two $\m$-primary ideals  $I_2 \subseteq I_1$ in a Cohen-Macaulay local ring of dimension $d \geq 1$ 
 (Theorem~\ref{thm-fundamental}). 
As a consequence we are able to describe the Hilbert coefficients   
 of the Hilbert polynomial associated to the function $\ell(R/I_{1} I_{2}^{n})$ which we denote by  $g_{i, I_1}(I_2)$ ($0 \leq i \leq d$)
(Lemma~\ref{hilb-coef-one}). 
A lower bound for $g_{1, I_1}(I_2)$ was given in \cite[Proposition~4.1]{jay-verma} under some assumptions on  
$I_2 \subseteq I_1$. In this paper,  we  improve their bound. We also give an upper bound for   $g_{1, I_1}(I_2)$.
  As a consequence,  we  show that when the lower bound is
attained $G_{I_1}(I_2)$ is 
Cohen-Macaulay  and when the upper bound is attained 
 $\depth~G_{I_1}(I_2) \geq d-1$ 
(Proposition~\ref{prop-al-min-mult-ineq}).

We  describe the 
 Hilbert coefficients  of the fiber cone $F_{I_1}(I_2)$,  
 $f_{i, I_1}(I_2)$, $0 \leq i \leq d-1$ 
(Lemma~\ref{cor-fiber-coeff}). 
The multiplicity of the fiber cone is of interest. 
An upper bound for the multiplicity  was given in
 \cite{cpv} and  in \cite{jay-verma} for $I_1=\m$. Using the complex 
$D(\x_{d}, {\Fi}, (1, n))$  ($n=1$)   obtain an upper  bound for the multiplicity in a more general setting  (Corollary~\ref{upper-bound-multiplicity}). We also give a lower bound on the multiplicity of the fiber cone.  In \cite{cpv} the authors remark that when the upper bound is attained the fiber cone need not be Cohen-Macaulay. We make an interesting observation. We show that when the upper bound is attained then $\depth ( (\x_d), G_{I_1}(I_2))=d$   and when the lower bound  is attained, $\depth ((\x_d), G_{I_1}(I_2)) \geq d-1$ (Corollary~\ref{upper-bound-multiplicity}). 
 
One interesting question is: How is the depth of the fiber cone and the associated graded ring related? We give an answer to this  question using homological methods. 
For ideals of minimal multiplicity and almost minimal multiplicity most of the homologies $H_i(C(\x_{d}{\mathcal F}, (1, n)))=0$ vanish giving a nice relation between the depth of $G(I_2)$ and $F_{I_1}(I_2)$ (Theorem~\ref{thm-jay-verma-almm}).  

We now describe the organization of this paper. 
In Section~\ref{complex},  we define the  $\C1n$ and  $\D1n$. We describe some interesting properties of these complexes. In Section~\ref{def-mm} we define ideals with minimal multiplicity and ideals with almost minimal multiplicity in terms of the homologies of $\C1n$.  In Section~\ref{vanishing-homologies} we compare the depths of the graded rings $F_{I_1}(I_2)$, $G_{I_1}(I_2)$ and $G(I_2)$.
In Section~\ref{hilbert-coefficients} we describe the Hilbert coefficients of $\ell(R/I_1 I_2^n)$. We also state a more general form of the Huneke's Fundamental Lemma. In Section~\ref{hilb-coef-fiber-cone} we describe the Hilbert coefficients of the fiber cone.
   In Section~\ref{hilb-fib-cone-mm} we describe Hilbert coefficients of $\ell(R/I_1 I_2^n)$ and of the fiber cone
 for ideals 
of   minimal multiplicity and ideals of almost minimal 
multiplicity.


\section{ The complex 
$\C1n$, $\D1n$ and Hilbert function }
\label{complex}

For any two ideals  $I_1$ and $I_2$ of  $R$, let 
$\F2i$ 
(resp. ${\Fi}$) be the $I_2$-filtration 
 $ \{ I_{2}^{n} \}_{n \in \Z}$
(resp. $\{I_{1} I_{2}^{n} \}_{n \in \Z}$).   We use the following convention: For any ideal $I$, $I^n=R$ if $n \leq 0$.
Let $k \geq 1$ and 
$\x_k:= x_1, \ldots, x_k$ be a sequence of elements in $I_2$. 
Using the   mapping cone construction, Marley and Huckaba constructed the following complex \cite{tom-huc}:
\begin{equation}
\label{newdisplay-1}
\Cnn:     
 0
\rar \f{R}{I_2^{n-k} }
\rar \cdots 
\rar \left( \f{R}{I_2^{n-i}} \right)^{k \choose i}
\rar \cdots 
\rar \left( \f{R}{I_2^{n-1}} \right)^{k \choose 1}
\rar \f{R}{I_2^{n}} \rar 0.
\end{equation}
In a   similar  way,  for the filtration $\Fi$,  we get the   following complex:
\begin{equation}
\label{newdisplay-2}
\C1n:
     0
\rar \f{R}{I_1 I_2^{n-k} }
\rar \cdots 
\rar \left( \f{R}{I_1 I_2^{n-i}} \right)^{k \choose i}
\rar \cdots 
\rar \left( \f{R}{I_1 I_2^{n-1}} \right)^{k \choose 1}
\rar \f{R}{I_1 I_2^{n}} \rar 0.
\end{equation}
The maps in (\ref{newdisplay-1}) and (\ref{newdisplay-2}) are induced by the Koszul complex $K.(\x_k;R)$.
Corresponding to the short exact sequence (\ref{newdisplay-2})  we  have the short  exact sequence of complexes:
{\mptmcomx
\beqn
     0 
\rar {{\Cnk1}}
\rar {{\C1n}}
\rar {{\Cn1k1}}[-1]
\rar 0
\eeqn}
and the corresponding long exact sequence of homologies:
{\mptmcomx
\beq
\label{main-homology} \nno
          \cdots
\rar &  H_i( \C1n )
\rar &    H_{i-1}(\Cn1k1)
\\ 
\rar      H_{i-1}(\Cnk1)
\rar&    H_{i-1}( \C1n )
\rar& \cdots .
\eeq
\mptm

\blem
\label{homology-four}
Let $k \geq 2$ and $n < k$. If  $(\x_k)  \subseteq I_1$, then 
 for all $  i -1 \geq n$,
\beqn
{ \displaystyle
H_{i}( {\C1n  } )
= \left( \f{R}{I_1} \right)^{k \choose i}
}.
\eeqn
\elem
\proof If $n<k$ and  $(\x_k) \subseteq I_1$, then  for all $i \geq n+1$ all the maps in the complex (\ref{newdisplay-2}) are zero. 
\qed

Using Lemma~\ref{homology-four}, when $(\x_k) \subseteq I_1$, we can  truncate  the complex $\C1n$  to get the complex
$\D1n$ where:
\begin{equation}
\label{newdisplay-3}
\begin{array}{rcll}
\D1n :&
     0
\rar \left( \f{R}{I_1 } \right)^{k \choose n}
\rar \cdots 
\rar \left( \f{R}{I_1 I_2^{n-i}} \right)^{k \choose i}
\rar \cdots 
\rar \left( \f{R}{I_1 I_2^{n-1}} \right)^{k \choose 1}
\rar \f{R}{I_1 I_2^{n}} \rar 0,
& \hspace{.1in}  n < k; \\
&
     0
\rar \f{R}{I_1 I_2^{n- k} }
\rar \cdots 
\rar \left( \f{R}{I_1 I_2^{n-i}} \right)^{k \choose i}
\rar \cdots 
\rar \left( \f{R}{I_1 I_2^{n-1}} \right)^{k \choose 1}
\rar \f{R}{I_1 I_2^{n}} \rar 0,
& \hspace{.1in} n \geq k;
\end{array} 
\end{equation}
By our construction, 
$H_i(\D1n)= 0$ for  $i >n$. In this section, we will show that  the complex  (\ref{newdisplay-3})  and the  complex
(\ref{newdisplay-1}) do share some similar properties. 

We  recall some basic facts and results. 
We say that an element 
$x \in I_2$  is superficial for 
$I_2$ and 
$I_1$ if there exists a positive integer 
$r_0$ such that for all 
$r \geq r_0$ and all $ s \geq 0$, 
$xR \cap I_1^s I_2^r  =x I_1^sI_2^{r-1}$. 
Let $k \geq 2$.  We say that
$\x_k$
 is a superficial sequence  if for each $i=1, \ldots,k$,
 $\overline{x_i} \in \overline{I_2}$ is superficial for $\overline{I_2}$ and $\overline{I_1}$ where 
 $\bar{\hphantom{sp}}$ denotes the image in $R/(x_1, \ldots, x_{i-1})$.
  Rees also showed that we can choose a minimal reduction
  $(\x_{d}) \in I_2$  which is  superficial sequence for $I_2$ and $I_1$  \cite{rees}. Hence, we can assume that $(\x_k)$ is a minimal reduction of $I_2$  and is generated by a superficial sequence for $I_2$ and $I_1$.

As the homologies of $\C1n$ and $\D1n$ can be computed using  techniques similar to those in  \cite{anna}, with slight modification,  we state the theorem  without proof.
 
\bt
\label{homology}
Let   $ k \geq 1$ and let $\x_{k} \in I_2$ be a superficial sequence for $I_2$ and $I_1$. Then
\been
\item
\label{homology-one}
For all $n \in \Z$, 
${\displaystyle
        H_0( \C1n    )
   \cong \f{R}
           {I_1I_2^{n} + (\x_{k})}.
}$\\
Let $(\x_d) \subseteq I_1$. Then ${\displaystyle
H_0(\D1n) = 
\left\{
\begin{array}{ll}
H_0(\C1n) & \hspace{.1in} \mbox{ for } n \geq 0,\\ 
0 & \hspace{.1in} \mbox{ for } n < 0.\\
\end{array} \right.
}$.

\item
\label{homology-two}
Let $\x_k$ be a regular sequence. Then for all $n \geq 2$, 
${\displaystyle
    H_1( {{\C1n                  }} )
\cong \f{(\x_{k}) \cap  I_1 I_2^{n}}
  {  (\x_{k}) I_1 I_2^{n-1}}}.
$\\
Let $(\x_d) \subseteq I_1$. Then ${\displaystyle
H_1(\D1n) = 
\left\{
\begin{array}{ll}
\f{(\x_{k}) \cap  I_1 I_2^{n}}
  {  \x_{k} I_1 I_2^{n-1}} & \hspace{.1in} \mbox{ for } n \geq 1,\\ 
0 & \hspace{.1in} \mbox{ for } n < 1.\\
\end{array} \right.
}.$

 \item
\label{homology-three}
For all $n \in \Z$, 
${\displaystyle
     H_{k}( {{\C1n  }} )
\cong \f{  I_1 I_2^{n-k+1} : ( \x_{k}) }
        {I_1 I_2^{n-k}}}.$\\
       Let $(\x_k) \subseteq I_1$. Then        ${\displaystyle
H_k(\D1n) = 
\left\{
\begin{array}{ll}
H_k(\C1n) & \hspace{.1in} \mbox{ for } n \geq k,\\ 
0 & \hspace{.1in} \mbox{ for } n < k\\
\end{array} \right.
}$.
\eeen
\et


The depth of the graded ring $G_{I_1}(I_2)$ and the vanishing of the homologies of the complex $\D1n$ are related as follows:
\bp
\label{grade}
Let $I_2 \subseteq I_1$ be $\m$-primary ideals in $(R, \m)$. Assume that $(\x_k)$ is a superficial sequence for $I_2$ and $I_1$ and that $(\x_k) \subseteq I_1$.  
For the filtration $\Fi$, 
\beqn 
\depth(( \x{_k})^{\star}, G_{I_1}(I_2)) = \min \{ H_{j-k}(\D1n) \not = 0  
\mbox { for some } n\}.
\eeqn
\ep 
\proof The proof is similar to the proof of Proposition~3.3
of \cite{tom-huc}. \qed

In the next lemma we show that the complex
$\D1n$ satisfies a certain rigidity similar to  that of the complex $\Cnn$.

\blem
\label{vanishing-rigidity}
Let $I_1, I_2$ and $\x_k$ be as in Proposition~\ref{grade}.
If $H_j(\D1n) = 0$  for some $j \geq 1$ and for all $n$, then 
$H_i(\D1n) = 0
$  for all $i \geq j$ and for all $n$.
\elem
\proof The proof follows by induction on $k$. \qed

We  state a crucial property satisfied by $\C1n$ and $\D1n$.

\blem
\label{vanishing}
Let $1 \leq i \leq k$. Let $\x_k \in I_2$ be a regular sequence in $R$ which is superficial for $I_2$ and $I_1$. Then for all $n \gg 0$ we have
$H_i(\C1n)=H_i(\D1n) = 0$. 
\elem
\proof Note that for $n \gg 0$, $H_i(\C1n)  = H_i(\D1n)$.
The proof follows by applying induction on $k$ for the complex $\D1n$. \qed

\bnot
   \beqn
   \begin{array}{lllllll}
     h_i ( C_{\bullet},\x_{k})(1,n) 
&:=& \ell
    \left( H_i (\C1n) 
    \right).
    &&
    h_i (C_{\bullet}, \x_{k})(0,n) 
&:=& \ell  \left( H_i (\Cnn)  \right).\\
 h_i ( D_{\bullet}, \x_{k})(1,n) 
&:=& \ell
    \left( H_i (\D1n) 
    \right).
    &&
    h_i (D_{\bullet},\x_{k}) (1, *)
&:=&  \sum_{n \geq 0} h_i (D_{\bullet},\x_{k}) (1, n).\\
 h_i (C_{\bullet},\x_{k}) (1, *)
&:=&  \sum_{n \in  \Z} h_i (C_{\bullet},\x_{k}) (1, n).
&&  h_i (D_{\bullet},\x_{k}) (1, *)
&:=&  \sum_{n \geq 0} h_i (D_{\bullet},\x_{k}) (1, n).
\end{array}
\eeqn
\enot

\bt
   \label{rigidity}
Let $(R, \m)$ be a  local ring. Let $I_1$ and $I_2$ be $\m$-primary ideals of $R$. Let  $k \geq 1$ and
$\x_{k}$ be a regular sequence which is superficial for $I_2$ and $I_1$.      
Let $i \geq 1$. If $(\x_k) \subseteq I_1$, then   
    \beq
    \label{rigidity-two}
 \sum_{j \geq i} (-1)^{j-i} h_j(D_{\bullet},\x_{k})(1, *) \geq 0.
\eeq
Equality holds if and only if
\beqn
\depth((\x_k)^{\star}, G_{I_1}(I_2)) \geq k-i+1.
\eeqn 
\et
\proof The proof is similar to the proof of Theorem~3.7
of \cite{tom-huc}. \qed


\bnot
\beq
\begin{array}{lllllll}
     \label{notation-hilbert}
      H_{\Fi}(1,n)
&:=& \ell \left( \f{R}{I_1 I_2^n} \right)
&\hphantom{sp}&
      H_{\Fi}^{\prime}(1,n)&:=& 
      \left\{
\begin{array}{ll} 
    \ell 
    \left( \f{R}{I_1 I_2^n} 
    \right) 
&   \mbox{ if } n \geq 0\\
     0 
& \mbox{ if } n <0.
\end{array}
\right. 
\end{array}
\eeq
\enot

We now give a relation between the homology modules and the Hilbert function $H_{\Fi}(1,n)$.

  \blem
\label{lemma-extension-fiber}
   Let $(R, \m)$ be a Cohen-Macaulay ring of dimension $d \geq 2$. 
   Let $I_1$ and $I_2$ be an $\m$-primary ideals of $R$ and 
    $(\x_{d}) \in I_2$ a minimal reduction which is a superficial  sequence for $I_2$ and $I_1$.
\been
\item
 For all $n \geq 2$,
\beq
       \label{alt-sum-new-1}   
 \sum_{i=0}^{d} (-1)^i{d\choose i} 
 H_{\Fi}(1,n-i)  =   e_0(  I_2)  
-     \ell 
       \left(  
       \f{I_1 I_2^{n}}
           {(\x_{d}) I_1 I_2^{n-1}} 
        \right)
        +\sum_{i=2}^{d} (-1)^i h_i(C_{\bullet},\x_{d})(1, n).
 \eeq

 \item
 Let $(\x_d) \subseteq I_1$. The  for all $n \in \Z$,
\beq
       \label{alt-sum-new-d}   
 \sum_{i=0}^{d} (-1)^i{d\choose i} 
 H_{\Fi}^{\prime}(1,n-i)  =   e_0(  I_2)  
-     \ell 
       \left(  
       \f{I_1 I_2^{n}}
           {(\x_{d}) I_1 I_2^{n-1}} 
        \right)
        +\sum_{i=2}^{d} (-1)^i h_i(D_{\bullet},\x_{d})(1, n).
 \eeq
\eeen
 \elem
\proof 
From Theorem~\ref{homology}, for all $n \geq 2$, 
\beq
\label{alt-sum-1} \nno
     h_0(C_{\bullet}, \x_{d})(1, n)
-    h_1(C_{\bullet}, \x_{d})(1, n)
&=& \ell \left(
    \f{R}{I_1 I_2^{n} + (\x_{d})} 
    \right)
- \ell \left( 
\f{( \x_{d}) \cap  I_1 I_2^{n}}
  {(\x_{d}) I_1 I_2^{n-1}}
  \right)\\  \nno
&=& \ell \left(
    \f{R}{ ( \x_{d})} 
    \right)
-   \ell \left( \f{I_1 I_2^{n} + (\x_{d})} 
                {( \x_{d})}
      \right)            
-     \ell \left( 
      \f{( \x_{d}) \cap I_2^{n}}
        {(\x_{d}) I_1 I_2^{n-1}}
       \right)\\
&=& e_0(I_{2}) 
-    \ell 
     \left( 
     \f{ I_1 I_2^{n}}
            { (\x_{d}) I_1   I_2^{n-1}}
     \right).
     \eeq
     For the complex $\D1n$, (\ref{alt-sum-1}) holds true for all $n \in \Z$.
\qed

 Lemma~\ref{homology-negative} gives us insight for the behaviour of homology modules for the complex $\C1n$ for $n \leq 0$. 
\blem
\label{homology-negative}
Let $I_1$ and $I_2$ be $\m$-primary ideals in a Cohen-Macaulay local ring $R$. Let $\x_k$ be a superficial sequence in $R$. 
\been
\item
\label{homology-negative-one}
  For $n \leq 0$, 
 $
 {\displaystyle
   \sum_{i=0}^{k} (-1)^i {k \choose i}
  h_i (C_{\bullet},\x_k)(1, n-i)=0.
  }$
  
  \item
  \label{homology-negative-two}
  We have a surjective map:
  \beq
   \label{homology-negative-two-eqn}
 \phi: \left( \f{R}{I_1} \right)^k 
  \lrar \f{(\x_k)}{I_1 (\x_k)}.
  \eeq
 
\item
\label{homology-negative-three}
Let $2 \geq i \geq k$. Then 
\beqn
\sum_{i=2}^k
(-1)^i
\left[ h_i (C_{\bullet},\x_k)(1, 1)
- {k \choose i} \ell \left( \f{R}{I_1} \right)
\right] \leq  0
\eeqn
Let $k=d$, then  equality holds if and only if $(\x_d) \subseteq I_1 $.
Moreover if equality holds, then the map $\phi$ in 
(\ref{homology-negative-two})
is an isomorphism.
\eeen  
  \elem
\proof For 
$n \leq 0$,   the complex (\ref{newdisplay-2}) 
is of the form
\beqn
    0
\rar \f{R}{I_1  }
\rar \cdots 
\rar \left( \f{R}{I_1 } \right)^{k \choose i}
\rar \cdots 
\rar \left( \f{R}{I_1 } \right)^{k \choose 1}
\rar \f{R}{I_1 } \rar 0.
\eeqn
Hence
    \beqn
   \sum_{i=0}^{k} (-1)^i {k \choose i}
  h_i(C_{\bullet}, \x_k) (1, n-i)
  = \sum_{i=0}^{k} (-1)^i {k \choose i}
 \ell \left( \f{R}{I_1} \right)=0.
  \eeqn
This proves(\ref{homology-negative-one}).

(\ref{homology-negative-two}) was proved in \cite[Lemma~3.4]{jay-verma-almm}.

We prove  (\ref{homology-negative-three}).
   The first part of (\ref{homology-negative-three}) follows from the fact that 
 $
 {\displaystyle 
 h_i (C_{\bullet},\x_k)(1, 1)
\leq {k \choose i} \ell \left( \f{R}{I_1} \right)
}.
$
We have 
\beq
\label{homology-negative-four}
&&  
      \ell \left( \f{R}      {I_1 I_2} \right)
-      \ell \left( \f{\x_d}   {I_1 (\x_d)} \right)
+      \sum_{i=2}^{d}  (-1)^i {d \choose i}
            \left( \f{R}      {I_1} \right)\\ \nno
&\geq& \ell \left( \f{R}      {I_1I_2} \right)
-       d~  \left( \f{R}      {I_1 } \right)
+      \sum_{i=2}^d (-1)^i {d \choose i}
            \left( \f{R}      {I_1} \right)
            \hspace{.5in} \mbox{by (\ref{homology-negative-two})}\\  \nno
&=&     e(\x_d) 
-      \ell \left( \f{I_1 I_2}{I_1 (\x_d)} \right)
+      \sum_{i=2}^d  (-1)^i 
        h_i(C_{\bullet}, \x_d) (1, 1) 
\hspace{.5in} \mbox{[putting $n=1$ in 
Theorem~\ref{lemma-extension-fiber}(1)].}
\eeq  
If $(\x_d)\subseteq I_1$, 
then  equality holds by Lemma~\ref{homology-four}. 
Conversely, suppose equality holds in 
(\ref{homology-negative-four}), then 
\beqn
h_i(C_{\bullet}, \x_d)(1,1) = {d \choose i} \ell \left( \f{R}{I_1} \right)
\hspace{.2in}
\mbox{for all }i=2, \ldots ,d.
\eeqn
In particular,
\beqn
\f{I_1 : (\x_d)}{I_1} = 
h_d(1,1) =  \ell \left( \f{R}{I_1} \right)
\eeqn
This implies that $I_1 : (\x_d) = R$. 

Finally,  if equality holds, in 
(\ref{homology-negative-four}), then  
$d~\ell \left( \f{R}{I_1} \right)
 = \ell \left( \f{(\x_d)}{I_1 (\x_d)} \right)$ 
 and hence the map in (\ref{homology-negative-two}) is an isomorphism.
\qed
\brm
If  $(\x_d) \not \subseteq I_1$, then 
$
\sum_{i=2}^d
(-1)^i
\left[ h_i (C_{\bullet}, \x_d)(1, 1) - \ell (R/I_1) \right]
$ may be non-zero (Example~\ref{countex}). 
\erm

 \bex
   \label{countex}
   Let $R= k[x,y,z]_{\m}$,  where $\m = (x,y,z)$. Let
 $I_2  = (x^3, y^3,z^3, xy, xz, yz)$
 and $(\x_3)= (x^3 + yz, y^3+ z^3+ xz, xz+ xy)$.
 Let  $n \geq 4$ and $I_1  = (x^n, y^3,z^3, xy, xz, yz)$.
 Here $(\x_3)$ is generated by a superficial sequence and $(\x_3)\not \subseteq I_1$.
The map in  (\ref{homology-negative-two-eqn}) is not an isomorphism.
\qed
\eex


\section{Minimal and almost minimal multiplicity and vanishing results}
\label{def-mm}

Goto \cite{goto} defined ideals of minimal multiplicity. Following Goto, Jayanthan and Verma defined ideals of minimal  (resp. almost minimal) multiplicity in \cite{jay-verma} (resp. \cite{jay-verma-almm})
 when $I_2 \subseteq I_1$. 
In this paper, we  generalize these definitions  for any 
two $\m$-primary ideals $I_2$ and $I_1$. 
We also show that, under some mild assumptions,   the  homologies of the complex $\C1n$ for these ideals  have nice vanishing properties.

\bdefn
\label{defn-mm-amm}
Let $(R, \m)$ be a Cohen-Macaulay local ring and let $I_1$ and  $I_2$   
 be $\m$-primary ideals in $R$. $I_2$ has  minimal multiplicity (resp. almost minimal multiplicity) with respect to $I_1$ if for some  minimal reduction $(\x_d)$ of 
$I_2$
\beqn
 \ell \left( \f{I_1 I_2}{I_1 (\x_d)} \right) 
 = 0
\hphantom{sp}\left( resp. \hphantom{sp} \ell \left( \f{I_1 I_2}{I_1 (\x_d)}  \right) 
 = 1 \right). 
\eeqn
   \edefn

\blem
\label{min-mult-defn}
Let $(R, \m)$ be a Cohen-Macaulay local ring of positive dimension. Then $I_2$ has minimal multiplicity (resp. almost minimal multiplicity) with respect to $I_1$ if and only if  for some minimal reduction $(\x_d)$ of $I_2$,
\beqn
e_0(I_2)
- \ell \left( \f{R}{I_1I_2} \right) 
+ d \cdot \ell \left( \f{R}{I_1} \right)
&=&  - \sum_{i=2}^d
(-1)^i
\left[ h_i(C_{\bullet}, \x_d) (1, 1)
- {d \choose i} \ell \left( \f{R}{I_1} \right)
\right]\\
\left( 
\mbox{resp. }e_0(I_2)
- \ell \left( \f{R}{I_1I_2} \right) 
+ d \cdot \ell \left( \f{R}{I_1} \right)
\right.
&=& 
\left.
1 - \sum_{i=2}^d
(-1)^i
\left[ h_i(C_{\bullet}, \x_d)(1, 1)
- {d \choose i} \ell \left( \f{R}{I_1} \right)
\right]\right).
\eeqn
\elem
 \proof We have:  
\begin{equation}
\label{newdisplay-4}
\Cd11:
      0
\rar \left( \f{R}{I_1} \right)^{d \choose d}
\rar \left( \f{R}{I_1} \right)^{d \choose d-1}
\rar \cdots
\rar \left( \f{R}{I_1} \right)^{d}
\rar \f{R}{I_1 I_2}
 \rar 0
\end{equation}
By Theorem~\ref{homology} we have
we have
\beqn
\ell \left( \f{R}{I_1 I_2} \right)
- d~\ell \left( \f{R}{I_1 } \right)
+ \sum_{i=2}^d
(-1)^i   {d \choose i} \ell \left( \f{R}{I_1} \right)
   &=& e_0(I_2)
- \ell \left( 
   \f{ I_1 I_2}
   {(\x_{d}) I_1}\right)
   + \sum_{i=2}^d(-1)^i   h_i (C_{\bullet}, \x_d)(1, 1).
\eeqn
Hence
\beqn
   e_0(I_2)
- \ell \left( \f{R}{I_1 I_2} \right)
+ d~\ell \left( \f{R}{I_1 } \right)
= - \sum_{i=2}^d(-1)^i  \left[
 h_i (C_{\bullet}, \x_d)(1, 1)
 - \ell \left( \f{R}{I_1} \right)
\right]
+ \ell \left( 
   \f{ I_1 I_2}
   {(\x_{d}) I_1}\right).
\eeqn
\qed

Example~\ref{countex-1} reveals that 
   $I_1 I_2= I_1 (\x_d)$ while
   $e_0(I_2) - \ell (I_1I_2)+ d~\ell (R/I_1) >1$.
Therefore,  the results in  \cite{jay-verma}  and 
\cite{jay-verma-almm} can be extended to ideals 
$I_2 \not \subseteq I_1$ if one takes into account the homologies 
$H_i(C_{\bullet}, \x_d) (1, 1)$, $i=2, \ldots, d$.

\bex
   \label{countex-1}
   Let $R= k[x,y,z]_{\m}$,  where $\m = (x,y,z)$. Let
 $I_2  = (x^2, y^2,z^2, xy, xz, yz)$
 and $(\x_3)= (x^2 + yz, y^2+ z^2+ xz, xz+ xy)$.
 Let  $n \geq 4$ and $I_1  = (x^n, y^2,z^2, xy, xz, yz)$.
 Here $(\x_3)$ is a minimal reduction of $I_2$, is generated by a superficial sequence and $(\x_3)\not \subseteq I_1$.
\qed
\eex

Using Lemma~\ref{homology-negative}
and  Lemma~\ref{min-mult-defn}  one can verify that if $(\x_d) \subseteq I_1$ for some minimal reduction of $I_2$,
  then $I_2$ has minimal (resp. almost minimal) multiplicity with respect to $I_1$  if and only if 
\beqn
e_0(I_2)
- \ell \left( \f{R}{I_1I_2} \right) 
+ d \cdot \ell \left( \f{R}{I_1} \right)
=  0\\
\left(resp.  \hphantom{sp} e_0(I_2)
- \ell \left( \f{R}{I_1I_2} \right) 
+ d \cdot \ell \left( \f{R}{I_1} \right)
=1\right).
\eeqn

The following lemma is useful as it  describes the vanishing of certain homology modules as well as describe the Hilbert coefficients of the fiber cone. 

\blem
\label{mm-amm}
Let $(R,\m)$ be a Cohen-Macaulay  local ring of dimension $d$. Let $I_1$ and $I_2$ be $\m$-primary ideals of $R$. Let $(\x_{d})$ be a minimal reduction of $I_{2}$.
\been
\item 
\label{lemma-mm-van}
Suppose $I_1I_2 = I_1 (\x_d)$. Then for all $n \geq 1$
\been
\item
\label{mm-1}
$
{\displaystyle
\ell \left( \f{I_1 I_2^n}{(\x_{d})I_1 I_2^{n-1}} \right)=0
}
$.

\item
\label{mm-3}
 $
  (\x_{d}) \cap I_1 I_2^n
=   (\x_{d}) I_1 I_2^{n-1}.
$
\eeen

\item
\label{lemma-amm-van}
Suppose  
$
{\displaystyle
\ell \left( \f{I_1 I_2}{(\x_{d})I_1 } \right)
=  1
}
$

\been
\item
\label{amm-1}
$
{\displaystyle
     \ell 
     \left( 
     \f{I_1 I_2^n}{(\x_{d})I_1 I_2^{n-1}} 
     \right)
\leq  1
}
$
for all $n \geq 1$.

\item
\label{amm-2}
For all $n \geq 1$, 
\beqn
     (\x_{d}) \cap I_1 I_2^n
&=& \left\{ 
\begin{array}{ll}
     (\x_{d}) I_1 I_2^{n-1} 
+    (y_1 y_2^n)  
&    \mbox{ if }  (\x_{d}): (y_1 y_2^n) = (1)\\
     (\x_{d}) I_1 I_2^{n-1} 
&   \mbox{ if } ( \x_{d}): (y_1 y_2^n) \subseteq \m
\end{array} \right. 
\eeqn
for some $y_{1}\in I_{1}$ and $y_{2} \in I_{2}$.
\eeen
\eeen
\elem
\proof (\ref{mm-1}) follows from the definition. Using 
 (\ref{mm-1}) we get
 $
\label{al-mm-2}
   (\x_{d}) \cap I_1 I_2^n
=  (\x_{d}) \cap(\x_{d}) I_1 I_2^{n-1}
=  (\x_{d}) I_1 I_2^{n-1}.
$
for all $n \geq 1$. 

\noindent
To prove (\ref{amm-1}) it is enough to show that for all $n\geq 1$
\beq
 \label{star-two}
   I_2I_2^n 
= (\x_{d}) I_1 I_2^{n-1} 
+ (y_1y_2^n), 
 \hspace{.2in}
 \mbox{and }  
 \hspace{.2in}
  \m   y_1 y_2 
\subseteq   (\x_{d}) I_1 I_2^{n-1}.
 \eeq
 Since
$
\ell (I_1I_2/ (\x_{d})I_1)
= 1$, 
there exists $y_1 \in I_1$ and $y_2 \in I_2$ such 
that 
\beq 
\label{star-one}
              I_1 I_2
&=&         ( \x_{d})I_1 + (y_1 y_2)\hspace{.2in} \mbox{ and } \hspace{.2in}         
    \m   y_1 y_2 
\subseteq   (\x_{d}) I_1.
\eeq
 If $n = 1$, then we are done by our assumption.
   Let  $n >1$. Clearly 
$(\x_{d}) I_1 I_2^{n-1} + (y_1 y_2^{n}) \subseteq  I_1 I_2^n$
  then by induction hypothesis
\beq \nno
   I_1 I_2^n  
&=&         (I_1 I_2^{n-1}) I_2 
=           (\x_{d}) I_1 I_2^{n-2} + (y_1 y_2^{n-1}) I_2
=           (\x_{d}) I_1 I_2^{n-1} +  (y_1 y_2^{n-1})I_2\\ 
\nno
            (y_1 y_2^{n-1}) I_2
&\subseteq&  (y_2)  I_1 I_2^{n-1}
=            (y_2) ((\x_{d}) I_1I_2^{n-2}  
+            (y_1 y_2^{n-1}))
\subseteq    (\x_{d}) I_1 I_2^{n-1} + (y_1 y_2^{n})\\
 \label{anni}
            \m  (y_1 y_2^n) 
&= &          \m  (y_1 y_2)   (y_2^{n-1}) 
\subseteq ( (\x_{d}) I_1 )(y_2^{n-1}) 
\subseteq ( \x_{d})  I_1 I_2^{n-1}. 
\eeq
This proves (\ref{amm-1}). 

We now prove (\ref{amm-2}).
 From (\ref{amm-1}) we have
\beq
\label{proof-amm-2}
(\x_{d}) \cap I_1 I_2^n
&=&       (\x_{d}) 
     \cap ((\x_{d}) I_1 I_2^{n-1} 
     + (y_1 y_2^n))\\ \nno
&=& (\x_{d}) I_1 I_2^{n-1} 
+  (\x_{d}) \cap (y_1 y_2^n)\\ \nno
&=& (\x_{d} )I_1 I_2^{n-1} 
+   (y_1 y_2^n) 
 ((\x_{d}): (y_1 y_2^n))
\eeq
Now two cases arise. \\
{\ncent Case(i):}$ (\x_{d}): (y_1 y_2^n)
= (1)$. In this case 
$(\x_{d}) \cap I_1 I_2^n
= (\x_{d}) I_1 I_2^{n-1} 
+   (y_1 y_2^n) 
$.\\
{\ncent Case(ii):}
$ (\x_{d}): (y_1 y_2^n) \subseteq \m$. Then by 
(\ref{amm-1}), 
$(\x_{d}) \cap I_1 I_2^n
\subseteq (\x_{d} )I_1 I_2^{n-1} 
+  \m (y_1 y_2^n) 
\subseteq (\x_{d}) I_1 I_2^{n-1}
$ by (\ref{anni}).
 This proves (\ref{amm-2}). 
\qed

The homologies of the complex $\Cdnn$ and $\Ddnn$ for ideals of minimal multiplicity and almost minimal multiplicity do satisfy some nice interesting vanishing properties. We list them in this section. 

\bp
\label{prop-van-almm}
Let $(R,\m)$ be a Cohen-Macaulay local ring, $I_1$ and $I_2$-primary ideals of $R$.  Let $\x_{d} \in I_2$ be a superficial sequence for $I_2$ and $I_1$. Suppose $(\x_d) \subseteq I_1$.
\been
\item
\label{prop-van-almm-mm}
Let  $I_2$ be  an ideal of  minimal multiplicity with respect to $I_1$.
Then 
\beqn
h_i(D_{\bullet}, \x_{d})(1, *) = 0,
 \hspace{.5in} \mbox{ for all }i \geq 1.
\eeqn

\item
\label{lemma-vanishing-ammm-2}
Let  $I_2$ be  an ideal of almost minimal  multiplicity with respect to $I_1$. 
\been
\item
For all $n \geq 1$,    
$
{\displaystyle
h_1(D_{\bullet},\x_{d})(1, n) \leq 1. }
$

\item
\label{lemma-vanishing-ammm-b}
If 
$(\x_{d})  \cap I_1 I_2= (\x_d)I_1$, 
then   
 $h_i(D_{\bullet}, \x_{d})(1, *)=0$ for all
 $i \geq 2$. 
\eeen \eeen
\ep
\proof
It is enough to show that $h_1(D_{\bullet}, \x_{d})(1, *) = 0$
[Theorem~\ref{vanishing-rigidity}]. 
By Theorem~\ref{homology}(\ref{homology-two})  and Lemma~\ref{mm-amm}, 
 for all $n \geq 1$, 
\beqn
        h_1( D_{\bullet},\x_{d})(1, n) 
\cong \f{(\x_{d}) \cap I_1 I_2^n}
          { (\x_{d}) I_1 I_2^{n-1}}
          =  \f{(\x_{d}) \cap I_1 I_2^n}
          {  I_1 (\x_d)^{n}}
          = 0,
          \eeqn
          if $I_2$ is  an ideal of  minimal multiplicity with respect to $I_1$.
          
 If $I_2$ is an ideal of  minimal multiplicity with respect to $I_1$,         then by Theorem~\ref{homology}(\ref{homology-two}) and Lemma~\ref{mm-amm},  for all $n \geq 1$, 
\beqn
        h_1(D_{\bullet},\x_{d})(1, n) 
\cong  \f{(\x_{d}) \cap I_1 I_2^n}
          {(\x_{d}) I_1 I_2^{n-1}} 
\cong \f{( \x_{d}) I_1 I_2^{n-1} + ((\x_{d}) \cap (y_1 y_2^n))}
        {\x_{d} I_1 I_2^{n-1}}
  \eeqn
  for some $y_1\in I_1$ and $y_2 \in I_2$.
By (\ref{anni}), 
$
           \m (y_1 y_2^n  \cap \x_d) 
 \subseteq \m y_1 y_2^n   
 \in       (\x_d) I_1 I_2^{n-1}
 $. 

(\ref{lemma-vanishing-ammm-b}) follows from 
Theorem~\ref{homology}(\ref{homology-two}) and Lemma~\ref{mm-amm}(\ref{lemma-amm-van}), Proposition~\ref{grade} and \cite[Theorem~4.2.1]{rossi-valla-2}.
\qed  


\section{ Comparing depth of $F_{I_1}(I_2)$, $G_{I_1}(I_2)$ and $G(I_2)$}
\label{vanishing-homologies}

Throughout this section we will assume that 
$(\x_d)$ is a minimal reduction of $I_2$  and $(\x_d)\subseteq I_1$. 
For all $n \geq 0$ we have the 
 exact sequence of complexes:
\beq
\label{koszul-fiber}
      0 &
\rar  & K({\x_{d}^o}, \Fi)(n)
\rar \D1n 
\rar  \Cnn
\rar 0
\eeq
where  $K({\x_{d}}^o, \Fi)$ is the Koszul complex of the fiber cone 
$F_{I_1}(I_2)$ with respect to the sequence $\x_{d}^o = x_{1}^o, \ldots, x_{d}^o$. We have the corresponding long exact sequence of complexes:
\mptmcomx
\beq
\label{exact-depth-2}
\begin{array}{lrrl}
     & \cdots
\rar & H_{i+1}(D( \x_{d}, {\mathcal F}, (1, n)))
\rar & H_{i+1}(C( \x_{d}, {\mathcal F}, (0, n)))\\
\rar & H_{i}  (K( \x_{d}^o, F_{I_1}(I_2)))(n)
\rar & H_{i}  (\D1n)
\rar & H_{i}  (\Cnn)
\rar \cdots
\end{array}.
\eeq
\mptm

\bt
[Depth Lemma]
\label{depth-lemma}
Let $I_1$ and $I_2$ be ideals in a local ring $(R, \m)$.
If $\depth~G(I_2) < \depth~G_{I_1}(I_2)$, then 
$\depth~F_{I_1}(I_2) = \depth~G(I_2)+1$.
%
%
\et
\proof
Use (\ref{koszul-fiber}), Proposition~\ref{grade} and Proposition~3.1 of \cite{tom-huc}.
\qed

Using the results in this paper we generalize Theorem~4.4 of \cite{jay-verma-almm} and give a simple proof of Proposition ~5.4 of \cite{jay-verma-almm}.

\bnot
Let $s(\x_d):= \min\{  n | (\x_d)I_1 I_2^n = I_1 I_2^{n+1} \}$.
\enot

\bt
\label{thm-jay-verma-almm}
Let $(R,\m)$ be a Cohen-Macaulay local ring, $I_1$ an $I_2$-primary ideals of $R$. 
    Let $(\x_{d})$ be a minimal reduction of $I_2$ which is generated by a superficial sequence. 
    \been
\item
\label{thm-jay-verma-almm-a}
Let $I_{2}$ be an ideal of minimal multiplicity with respect to $I_{1}$. 
\been
\item
\label{thm-jay-verma-almm-a-i}
$G_{I_1}(I_2)$ is Cohen-Macaulay;

\item
\label{thm-jay-verma-almm-a-ii}
For all $i=0, \ldots, d-1$, 
$\depth( F_{I_1}(I_2))= i+1$ if and only if 
$\depth(  G(I_2)) =i$.
\eeen

\item
\label{thm-jay-verma-almm-b}
Let $I_{2}$ be an ideal of almost minimal multiplicity with respect to $I_{1}$.
 Suppose  $(\x_d) \cap I_1 I_2 = I_1 (\x_d)$. 
\been
\item
\label{thm-jay-verma-almm-b-i}
For all $i=1, \ldots, d-2$, $\depth((\x_{d}), F_{I_1}(I_2))= i+1$ if and only if 
$\depth((\x_{d}), G(I_2)) =i$.

\item
\label{thm-jay-verma-almm-c}
 $F_{I_{1}}(I_{2})$ is Cohen-Macaulay if and only if 
 $\depth((\x_{d}), G(I_2)) \geq d-1$
 and
$\ell(I_{1}I_{2}^{n} + (\x_d) I_{2}^{n-1}/ 
 (\x_d)I_{2}^{n-1})=1$
     for all
$n = 1,\ldots, s(\x_d)$.
\eeen
\eeen
\et
\proof 
 (\ref{thm-jay-verma-almm-a-i} )
follows from
Proposition~\ref{prop-van-almm}.
(\ref{thm-jay-verma-almm-a-ii}) follows form
(\ref{thm-jay-verma-almm-a-i}) and Theorem~\ref{depth-lemma}(1).

%

 (\ref{thm-jay-verma-almm-b-i}) 
follows form \cite[4.2.1]{rossi-valla-2}, 
Proposition~\ref{prop-van-almm} and Lemma~\ref{depth-lemma}(1) .

We now prove (\ref{thm-jay-verma-almm-c}).
  Applying Theorem~\ref{homology} and Lemma~3.2 in \cite{tom-huc} to the above sequence we get 
  \beqn
       \cdots
\rar H_2(C( \x_{d}, G(I_2)))(n)
 \rar  H_{1}  (K( \x_{d}^o, F_{I_1}(I_2)))(n)
  \rar  \f{(\x_d) \cap I_{1} I_{2}^{n}}
          {(\x_d)I_{1}I_{2}^{n-1}} 
  \rar  \f{(\x_d) \cap  I_{2}^{n}}
          {(\x_d)I_{2}^{n-1}} 
  \rar  \cdots.
  \eeqn
  Thus we have the exact sequence
 \beqn
       \cdots
\rar H_2(C( \x_{d}, G(I_2)))(n)
 \rar  H_{1}  (K( \x_{d}^o, F_{I_1}(I_2)))(n)
  \rar  \f{(\x_d) I_2^{n-1 }\cap I_{1} I_{2}^{n}}
          {(\x_d)I_{1}I_{2}^{n-1}} 
          \rar 0
          \eeqn 
          Hence $H_{1}  (K( \x_{d}^o, F_{I_1}(I_2)))(n)=0$ if and only if  $H_2(C( \x_{d}, G(I_2)))(n)$
and
${\displaystyle
  \ell \left(
     \f{(\x_d)I_{2}^{n-1} \cap I_{1} I_{2}^{n}}
       {(\x_d)I_{1}I_{2}^{n-1}} 
      \right)=0.}$
         Moreover for $n=1, \ldots, s(\x_d)$,
      \beqn
      1
 &=&   \ell \left(
     \f{I_{1} I_{2}^{n}}
        {(\x_d)I_{1}I_{2}^{n-1}} 
     \right)\\
&=&    \ell
     \left( 
     \f{I_{1} I_{2}^{n}}
       {(\x_d)I_{2}^{n-1} \cap I_{1} I_{2}^{n}}
     \right)
+    \ell \left(
     \f{(\x_d)I_{2}^{n-1} \cap I_{1} I_{2}^{n}}
       {(\x_d) I_{2}^{n-1}} 
     \right)\\
&=&      \ell 
     \left( 
     \f{I_{1} I_{2}^{n} + (\x_d)I_{2}^{n-1} }
       {(x_d) I_{2}^{n-1}}\right)
+    \ell \left(
     \f{(\x_d)I_{2}^{n-1} \cap I_{1} I_{2}^{n}}
       {(\x_d)I_{1}I_{2}^{n-1}} 
     \right)\\
&=&     \ell 
     \left( 
     \f{I_{1} I_{2}^{n} + (\x_d)I_{2}^{n-1} }
       {(\x_d) I_{2}^{n-1}}\right)
       +    \ell (H_{1}  (K( \x_{d}^o, F_{I_1}(I_2)))(n)  ) .
 \eeqn
 This gives 
$ \ell (H_{1}  (K( \x_{d}^o, F_{I_1}(I_2)))(n)  ) =0$
 if and only if 
 ${\displaystyle 
  \ell 
     \left( 
     \f{I_{1} I_{2}^{n} + (\x_d) I_{2}^{n-1} }
       {(\x_d) I_{2}^{n-1}}\right)=1}$.
\qed

\section{  Hilbert coefficients}
\label{hilbert-coefficients}
Throughout this section we will assume that  $(R, \m)$ is  a 
Cohen-Macaulay local ring. 
Let $I$ be an $\m$-primary ideal of $R$.  It is well known that for $n \gg 0$, the function $H(n):= \ell(R/I^n)$ is a polynomial in $n$ and  we will denote by $P(n)$. For a two dimensional Cohen-Macaulay local ring and an $\m$-primary ideal $I$, Huneke gave a relation between $\Delta^2[P(n) - H(n)]$ and the multiplicity of the ideal $I$ which is known as 
 Huneke's fundamental (\cite[Lemma~2.4]{huneke}). 
 This was generalized for an $\m$-primary ideal in a d-dimensional Cohen-Macaulay local ring in \cite{huc} and for an Hilbert filtration in \cite{tom-huc}.

 In \cite[Proposition~2.5]{jay-verma-almm}, Huneke's fundamental lemma was generalized for the filtration ${\mathcal F} = \{I_1I_2^n\}_{n \geq 0}$ for a two dimensional Cohen-Macaulay ring. We generalize Huneke's fundamental lemma for the filtration 
 $\Fi$  and for any dimension $d \geq 1$ (Theorem~\ref{thm-fundamental}). If we put $I = I_1 = I_2$ then we can recover Huneke's result as well as the result of Huckaba and Marley.

For all $n \in \Z$ let 
$H_{\Fi}(1,n):= \ell (R/ I_1 I_2^n)$ be the Hilbert function of $\F2i$ and let $P_{\Fi}(1,n)$ denote the corresponding 
Hilbert  polynomial. This polynomial can be written in the form
\beqn
\label{g-i-2}
   P_{\Fi}(1,n)
= \sum_{j=0}^{d} (-1)^j
   g_{j, I_1}(I_2) 
   {n + d-j -1 \choose d-j}.
\eeqn

 We now describe the coefficients $g_{i, I_1}(I_2)$, $0 \leq i \leq d$. Our result  is analogous  to \cite[Lemma~2.8, Proposition~2.9]{huc}.

\bp
\label{proposition-alt-sum-coeff}
      Let $I_1$ and $I_2$ be $\m$-primary ideals in a local ring $(R, \m)$ of dimension $d \geq 1$. Let $(\x_{d})$  be a minimal reduction which is generated by a superficial sequence for $I_2$ and $I_1$. Let  $1 \leq i \leq d$.
      \been
      \item
       \label{proposition-alt-sum-coeff-1}
$
{\displaystyle
      \Delta^{d-i} 
      \left[ P_{\Fi}(1,0) 
      \right]
=      (-1)^i g_{i, I_1}(I_2).
      }$
      
      \item
      \label{proposition-alt-sum-coeff-2}
${\displaystyle
    g_{i, I_1}(I_2) 
=     \sum_{n \geq  i-1}      {n \choose i-1}
   \Delta^{d} 
   \left[ P_{\Fi}(1,n+1) 
        - H_{\Fi}(1,n+1) \right] 
        + (-1)^d \ell \left( \f{R}{K_i} \right), }$
where $K_i= R$ if $i< d$ and $K_d = I_1$.
\eeen
\ep
\proof
Let $d \geq 1$. If $i=d$, then  we have 
$ \Delta^{d-i} P_{\Fi}(1,0)
= \Delta^0 P_{\Fi}(1,0) 
=   P_{\Fi}(1,0) 
=  (-1)^dg_{d,{I_1}}(I_2)$.
 Now let $i<d$ and let $\x \in I_2$ be superficial for $I_1$ and $I_2$.
  If $\overline{\hphantom{xx}}$ denotes the image in $R/x$, then 
\beqn
 \Delta^{d-i} 
         P_{\Fi}(1,0)
      =  \Delta^{d-i-1} 
          [P_{\oFi}(1,0)]
          = (-1)^i  g_{i, \overline{I_1}}(\overline{I_2}) 
          = (-1)^i  g_{i,{I_1}}({I_2}) .
  \eeqn
  This proves (\ref{proposition-alt-sum-coeff-1})

 We now prove (\ref{proposition-alt-sum-coeff-2}). 
 \beqn
&&    \sum_{n \geq  i-1}      {n \choose i-1}
      \Delta^{d} 
      \left[ P_{\Fi}(1,n+1) 
           - H_{\Fi}(1,n+1) \right] \\
&=&  (-1)^{i-1} 
       \sum_{n \geq 0}     
       \Delta^{d-i+1} 
       \left[ P_{\Fi}(1,n+1) 
            - H_{\Fi}(1,n+1) 
       \right]
       \mbox{ \cite[Lemma~2.7]{tom-huc}  }               \\
&=&     (-1)^{i} 
       \Delta^{d-i}  
       \left[ P_{\Fi}(1,0) 
            - H_{\Fi}(1,0) 
       \right] \\
&=&     g_{i, I_1}(I_2) 
+      (-1)^{i+1} 
      \Delta^{d-i} H_{\Fi}(1,0)
      \hspace{3.1in}
     [\mbox{ by } (\ref{proposition-alt-sum-coeff-1})].  
 \eeqn
One can verify that 
\beqn
    \Delta^{d-i} H_{\Fi}(1,0)
&=& \left\{ 
    \begin{array}{ll}
    \ell \left( \f{R}{I_1} \right) & \mbox { if } i=d\\
     0      
&   \mbox { if } i< d\\
    \end{array}
    \right.
\eeqn
\qed          
                  
\blem
\label{hilb-coef}
Let $(R, \m)$ be a Cohen-Macaulay 
ring of dimension $d \geq 1$. 
Let $I_1$ and $I_2$ be an $\m$-primary ideals of $R$ and 
$(\x_{d})$ be a minimal reduction of $I_2$ which is generated by a superficial sequence for $I_2$ and $I_1$. Then 
      $
\displaystyle
   {g_{0, I_1}(I_2) 
 =  e ( I_2).}$
 \elem
\proof From  Lemma~\ref{vanishing} we get 
 $h_i( \x_{d})(1, n) = 0$ for all $i \geq 1$ and for all $n \gg 0$. Since  $H_{\mathcal F}(1,n)$ is a polynomial $P_{\mathcal F}(1,n)$ for all $n \gg 0$, 
      from Lemma~\ref{lemma-extension-fiber}
   we have
 \beq
\label{coefficient-fiber-g1}
&& g_{0, I_1}(I_2) \\ \nno
&=& \limn
    \sum_{i=0}^{d-1} (-1)^i{d-1 \choose i} 
      P_{\Fi}(1,n-i)\\ \nno
&=& \limn
    \sum_{i=0}^{d-1} (-1)^i{d-1 \choose i} 
     H_{\Fi}(1,n-i)
=    e_0(I_2)  
-   \limn~\ell 
    \left(  
    \f{I_1 I_2^{n}}
      {(\x_{d}) I_1 I_2^{n-1}} 
    \right)
=    e_0( I_2).
\eeq
The last equality follows from the fact that $(\x_{d})$ is a minimal reduction of $I_2$ and hence  $I_1 I_2^{n}
       =(\x_{d}) I_1 I_2^{n-1}$ for all $n \gg 0$.  
This proves the lemma.
\qed

\blem
[Huneke's fundamental lemma for ${\Fi}$]
\label{thm-fundamental}
 Let $(R, \m)$ be a Cohen-Macaulay ring of dimension $d \geq 2$. 
   Let $I_1$ and $I_2$ be an $\m$-primary ideals of $R$ and 
    $(\x_{d})$ be a minimal reduction of $I_2$ which is generated by a superficial sequence for $I_2$ and $I_1$. For all 
    $n \geq 1$,
\beqn
\Delta^d \left[ P_{\Fi}(1,n) - H_{\Fi}(1,n) \right]
&=&  \ell 
     \left( 
     \f{ I_1 I_2^{n} + (\x_d)}
       { (\x_{d})}
     \right)   
-    \sum_{i = 1}^d (-1)^i 
      h_i(C_{\bullet}, \x_{d})(1, n) \\
&=&  \ell 
     \left( 
     \f{ I_1 I_2^{n}}
       { (\x_{d}) I_1   I_2^{n-1}}
     \right)   
-    \sum_{i = 2}^d (-1)^i 
      h_i(C_{\bullet}, \x_{d})(1, n) .
     \eeqn
\elem
\proof For all $n \geq 0$,  
$ \Delta^d \left[P_{\mathcal F}(1,n) \right] 
=  g_{0, I_1}(I_2)$. 
Now  use Lemma~\ref{hilb-coef}
and Lemma~\ref{lemma-extension-fiber}.
\qed

We are ready to  describe the coefficients 
$g_{i, I_1}(I_2)$ explicitly, in terms of the homology modules of the complex $\C1n$.
 
      \blem
      \label{hilb-coef-one}
     Let $(R, \m)$ be a Cohen-Macaulay ring of dimension   $d \geq 1$. Let $I_1$ and $I_2$ be $\m$-primary ideals
  and 
    $(\x_{d}) $  a minimal reduction of $I_2$   which is generated by a superficial sequence for $I_2$ and $I_1$. Let $1 \leq i \leq d$. We can write:
    \mptmcomx
    \beq
 \label{hilb-coef-one-eq}
    g_{i, I_1}(I_2)
&=& 
   \sum_{n \geq  i-1}      {n \choose i-1}      
   \left[ 
   \ell 
   \left( 
   \f{I_1 I_2^{n+1} +(\x_{d}) }
     {(\x_{d}) } 
   \right)
-  \sum_{j \geq 1} (-1)^j
    h_j( C_{\bullet}, \x_{d})(1, n+1)\right] 
+ (-1)^d\ell \left( \f{R}{K_i} \right)
\\ 
%
\label{hilb-coef-one-eq-two}
&=& 
  \sum_{n \geq  i-1}      {n \choose i-1}      
  \left[ 
  \ell 
  \left( 
  \f{I_1 I_2^{n+1}}
    {(\x_{d}) I_1 I_2^{n}} 
  \right)
- \sum_{j \geq 2} (-1)^i
   h_j( C_{\bullet}, \x_{d})(1, n+1)\right]
+ (-1)^d\ell \left( \f{R}{K_i} \right)
\eeq
where $K_i = R$ for $1 \leq i \leq d-1$ and $K_d= I_1$.
\elem 
\mptm
\proof
The proof follows from Proposition~\ref{proposition-alt-sum-coeff}(\ref{proposition-alt-sum-coeff-2}) and Lemma~\ref{thm-fundamental}.
\qed

A formula for $g_{1, I_1}(I_2)$ plays a very important role in analyzing the depth of $G_{I_1}(I_2)$. When $(\x_d) \subseteq I_1$, we have nice formulas in terms of the homology modules. 
%

\blem
\label{formula-g1}
Let $(\x_d)\subseteq I_1$. With the assumptions as in 
Lemma~\ref{hilb-coef-one}
\beq
\label{formula-g1-one}
g_{1, I_1}(I_2)
&=& \sum_{n \geq 1}
    \ell 
     \left( \f{I_1 I_2^{n} + (\x_d)}
              {(\x_{d}) } 
     \right)
-    \sum_{j \geq 1} (-1)^{j} h_j(D_{\bullet},\x_{d})(1, *)
-    \left( \f{R}{I_1} \right)\\
    \label{formula-g1-two}
&=& \sum_{n \geq 1}
    \ell 
     \left( \f{I_1 I_2^{n}}
              {(\x_{d}) I_1 I_2^{n-1}} 
     \right)
-    \sum_{j \geq 2} (-1)^{j} h_j(D_{\bullet},\x_{d})(1, *)
-    \left( \f{R}{I_1} 
     \right). 
     \eeq
\elem
\proof Put $i=1$ in equation (\ref{hilb-coef-one-eq}) of Lemma~\ref{hilb-coef-one}. We get
\beqn
g_{1, I_1}(I_2)
&=& \sum_{n \geq 1}      
    \left[ \ell 
    \left( 
    \f{I_1 I_2^{n} +(\x_{d}) }
       {(\x_{d}) }  \right)
-   \sum_{j \geq 1} (-1)^j
     h_j( C_{\bullet}, \x_{d})(1, n)\right]\\
&=& \sum_{n \geq 1}      
    \left[ \ell 
    \left( \f{I_1 I_2^{n} +(\x_{d})}{(\x_{d})}\right)
-   \sum_{j \geq 1} (-1)^j
    \left[ \sum_{n \geq j}
      h_j( C_{\bullet}, \x_{d})(1, n)
+   \sum_{n=1}^{j-1} h_j( C_{\bullet}, \x_{d})(1, n)
    \right] \right]
                \eeqn
For all $i \geq 2$,
\beqn
 h_j(C_{\bullet},\x_{d})(1, n)
 = \left\{
 \begin{array}{ll}
  h_j(D_{\bullet},\x_{d})(1, n) & \mbox{ if } n \geq j\\
  {d \choose j} \ell \left( \f{R}{I_1} \right) & \mbox{ if } 1 \leq n \leq j-1.
 \end{array}
 \right.
 \eeqn
      Therefore,
\beqn
     g_{1, I_1}(I_2)        
&=& \sum_{n \geq 1}
    \ell 
    \left( \f{I_1 I_2^{n} + (\x_d)}
             {(\x_{d})} 
    \right)
-   \sum_{j \geq 1} (-1)^{j} 
      h_j(D_{\bullet},\x_{d})(1, *)
-    \left[ 
     \sum_{j=2}^d (-1)^j (j-1) {d \choose j} 
     \left( \f{R}{I_1} 
     \right)
     \right]
                \\
 &=& \sum_{n \geq 1}\ell 
     \left( \f{I_1 I_2^{n} + (\x_d)}
              {(\x_{d})} 
     \right)
-    \sum_{j \geq 1} (-1)^{j} h_j(D_{\bullet},\x_{d})(1, *)
-    \left( \f{R}{I_1} 
     \right). 
\eeqn
In a similar way, we can prove (\ref{formula-g1-two}).
\qed

If $(\x_d)$ is a minimal reduction of $I_2$  and 
$(\x_d) \subseteq I_1$, then we  can give   bounds on 
 $g_{1, I_1}(I_2)$. In  \cite[Proposition~4.1]{jay-verma} only a lower bound was given. We also give an upper bound. 
\bp
\label{prop-al-min-mult-ineq}
 Let $(\x_d)\subseteq I_1$. With the assumptions as in 
Lemma~\ref{hilb-coef-one}:
\been
\item 
${\displaystyle
\label{al-min-mult-eq-one} 
 g_{1, I_1}(I_2)
 \geq
      \sum_{n \geq 1}
      \left( \f{I_1 I_2^n+ (\x_{d})}
               {(\x_{d}) }
      \right)  
-     \ell \left( \f{R}{I_1} \right) }
$
and equality holds if and only if $G_{I_1}(I_2)$ is Cohen-Macaulay.

\item
\label{al-min-mult-eq-two} 
${\displaystyle 
 g_{1, I_1}(I_2)
\leq  \sum_{n \geq 1}
      \left( 
      \f{I_1 I_2^n}
        {(\x_{d}) I_1 I_2^{n-1}} \right)
-     \ell \left( \f{R}{I_1} \right) }
$
and equality holds if and only if $\depth~G_{I_1}(I_2) \geq d-1$.
\eeen
\ep

\proof The proof follows from Theorem~\ref{rigidity},
Proposition~\ref{grade} and
Proposition~\ref{prop-al-min-mult-ineq}
\qed

\bco
Let $(\x_d)\subseteq I_1$. With the assumptions as in 
Lemma~\ref{hilb-coef-one}: 
\been
\item
$g_{1, I_1}(I_2) \geq - \ell ( R/I_1)$ and equality holds if and only if $G_{I_1}(I_2)$ is Cohen-Macaulay and $I_1 I_2^n  + (\x_d)= (\x_d)$ for all $ \geq 1$.
 
\item
${\displaystyle
      \ell 
      \left( 
      \f{R}{I_1I_2} \right) 
\geq   e_0(I_2) - g_{1, I_1}(I_2)
+     \ell 
      \left( 
      \f{I_1I_2+ (\x_d)}{I_1 I_2}
\right)}$
and equality holds if and only if $G_{I_1}(I_2)$ and 
$I_1 I_2^n+(\x_d) = (\x_d)$  for all $n \geq 2$.
\eeen
\eco
\proof Both (1) and (2) follows from
Proposition~\ref{prop-al-min-mult-ineq}(1).
                        
\section{ Hilbert coefficients of the fiber cone}
\label{hilb-coef-fiber-cone}
Throughout this section we will assume $I_1$ and $I_2$ are $\m$-primary ideals in a Cohen-Macaulay local ring $(R,\m)$In this section we describe the Hilbert coefficients of the fiber cone in terms of the length of the homologies of the complex $\C1n$ and $\Cnn$.

\bnot 
We denote the Hilbert function (resp. polynomial) of the fiber cone by
\beqn
H ( F_{I_{1}}(I_{2}), n) 
= \ell \left( \f{I_{2}^{n}}{I_{1} I_{2}^{n}}\right)
\hphantom{space}
\left( \mbox{resp. }
P ( F_{I_{1}}(I_{2}), n)  
= \sum_{i=0}^{d-1} (-1)^i 
f_{i, I_1}(I_2) 
{n + d-1-i \choose d-1-i} \right). 
\eeqn
\enot

\brm
\label{hilbert-fiber cone}
Since
\beqn
    \ell \left( \f{I_2^n}{I_1 I_2^n}\right)
&=& \ell \left( \f{R}{I_1 I_2^n}\right)
- \ell \left( \f{R}{I_2^n}\right)\\
&=& \sum_{i=0}^{d-1} (-1)^{i}
    \left[e_{i+1}(I_{2}) - g_{i+1, I_1}{(I_2)}
+   e_{i}(I_{2}) - g_{i, I_1}(I_2) \right]
    {n + d-1 -i\choose d-1-i},
\eeqn
 we have 
\beqn
  f_{i, I_1}(I_2) 
= e_{i+1}(I_{2}) - g_{i+1, I_1}(I_2) 
+  e_{i}(I_{2}) - g_{i, I_1}(I_2) \hspace{.2in} 
\mbox{ for all }
\hphantom{space} i=0, \ldots, d-1.
\eeqn
\erm

\blem
\label{cor-fiber-coeff}
     Let  $d \geq 1$ 
  and 
    $(\x_{d})$  a minimal reduction of $I_2$ which is generated by a superficial sequence for $I_2$ and $I_1$. Let $1 \leq i \leq d$. We can write:
\mptmcomx
\beqn
  f_{i, I_1}(I_2) 
  = \sum_{n \geq i} {n \choose i}
    \left[ 
    \ell 
    \left( 
    \f{I_2^{n}}
      {I_1 I_2^{n} +  (\x_d) \cap I_2^n }
    \right)
   -  \sum_{j \geq 1} (-1)^j \left[
    h_j( C_{\bullet}, \x_{d})(0, n)
    - h_j( C_{\bullet}, \x_{d})(1, n)
    \right]
    \right]
- (-1)^d \ell \left( \f{R}{K_{i}} \right)
    \eeqn
    \beqn
  f_{i, I_1}(I_2) 
  = \sum_{n \geq i} {n \choose i}
    \left[ \left[
    \ell 
    \left( 
    \f{I_2^{n}}
      {I_1 I_2^{n}}
    \right)
    -\ell 
    \left( 
    \f{(\x_d)I_2^{n-1}}
      {(\x_d)I_1 I_2^{n-1}}
    \right)
        \right]
   -  \sum_{j \geq 2} (-1)^j \left[
    h_j( C_{\bullet}, \x_{d})(0, n)
    - h_j( C_{\bullet}, \x_{d})(1, n)
    \right]
    \right]
- (-1)^d \ell \left( \f{R}{K_{i}} \right)
    \eeqn
    where $K_i = R$ for $i=0, \ldots, d-2$ and $K_{d-1} = I_1$
    \elem
    \proof The proof follows from (4.5) of \cite{tom-huc} and 
  Lemma~\ref{hilb-coef-one}.  

\mptm
We give  bounds for the multiplicity of the fiber cone. This is an improvement of \cite{cpv} and \cite[Corollary~4.2]{jay-verma}. 
\bco
\label{upper-bound-multiplicity}
Let $I_1$ and $I_2$ be $\m$-primary ideals in a local ring $(R, \m)$ of dimension at least two. Let $(\x_{d})$  be a minimal reduction which is a superficial sequence for $I_1$ and $I_2$. Suppose  
$(\x_d) \subseteq I_1$ and $(\x_d) \cap I_1 I_2 = (\x_d) I_1$. 
\been
\item
\label{upper-bound-0}
${\displaystyle
f_{0, I_1}(I_2) 
\leq  e_{1}(I_{2})  -  e_{0}(I_{2}) 
+ \ell \left( \f{I_2}{I_1 I_2} \right) 
+ \ell \left( \f{R}{I_2} \right)
-   (d-1)~\left( \f{R}{I_1} \right)} $
and equality holds if and only if 
$G_{I_1}(I_2)$ is Cohen-Macaulay and for all $n \geq 2$, 
 $I_1 I_2^n + (\x_d) = (\x_d)$.
 \vspace{.2in}
 
 \item
\label{upper-bound-1}
${\displaystyle
f_{0, I_1}(I_2) 
\geq  e_{1}(I_{2})
-      \left( 
       \f{I_1 I_2}
         {(\x_{d}) I_1 } 
       \right)
-      \sum_{n \geq 2}
       \left( 
       \f{I_1 I_2^n}
         {(\x_{d}) I_1 I_2^{n-1} } 
       \right)
+     \ell \left( 
      \f{R}{I_1} \right)} $
and equality holds  if and only if 
$\depth~G_{I_1}(I_2) \geq d-1$.
\eeen
\eco
\proof By Remark~\ref{hilbert-fiber cone} and Proposition~\ref{prop-al-min-mult-ineq} we get
\beqn
        f_{0, I_1}(I_2) 
&=&     e_{1}(I_{2}) - g_{1, I_1}(I_2)\\
&\leq&  e_{1}(I_{2})
-      \sum_{n \geq 1}
       \left( \f{I_1 I_2^n + (\x_d)}
                {(\x_{d}) } 
       \right)
  +    \ell \left( \f{R}{I_1} \right)\\
&\leq&     e_{1}(I_{2})
-      \left( 
       \f{I_1 I_2  +(\x_d)}
         {(\x_{d}) } 
       \right)
        +    \ell \left( \f{R}{I_1} \right)\\
        &\leq&     e_{1}(I_{2})
-      \left( 
       \f{I_1 I_2 }
         {(\x_{d} )I_1 } 
       \right)
        +    \ell \left( \f{R}{I_1} \right).
\eeqn
Using Lemma~\ref{homology-negative} we get:
\beqn
     \left( \f{I_1 I_2}{(\x_{d}) I_1 } \right)
=     e(\x_d) 
-    \left(  \f{R}{I_1 I_2} \right)
+  d~\left( \f{R}{ I_1}     \right)
=     e(\x_d) 
-    \left(  \f{I_2}{I_1 I_2} \right)
-    \left( \f{R}{I_2}      \right)
+  d~\left( \f{R}{ I_1}       \right). 
\eeqn
If equality holds if and only if  
${\displaystyle
     g_{1, I_1}(I_2)
=   \left( 
    \f{I_1 I_2}
      {(\x_{d}) I_1} 
    \right)
-   \ell 
    \left( 
    \f{R}{I_1} 
    \right)
}$ and
$I_1 I_2^n + (\x_d) = (\x_d)$ for all $n \geq 2$.
Now apply Proposition~\ref{prop-al-min-mult-ineq}(1).
This proves (\ref{upper-bound-0}).

Once again by Remark~\ref{hilbert-fiber cone} and Proposition~\ref{prop-al-min-mult-ineq} we get
\beqn
        f_{0, I_1}(I_2) 
&=&     e_{1}(I_{2}) - g_{1, I_1}(I_2)\\
&\geq&  e_{1}(I_{2})
-      \sum_{n \geq 1}
       \left( \f{I_1 I_2^n}
                {(\x_{d})I_1 I_2^{n-1} } 
       \right)
  +    \ell \left( \f{R}{I_1} \right)\\
&=&     e_{1}(I_{2})
-      \left( 
       \f{I_1 I_2}
         {(\x_{d}) I_1 } 
       \right)
-      \sum_{n \geq 2}
       \left( 
       \f{I_1 I_2^n}
         {(\x_{d}) I_1 I_2^{n-1} } 
       \right)
+     \ell \left( 
      \f{R}{I_1} \right).
 \eeqn
Applying 
  Proposition~\ref{prop-al-min-mult-ineq}(2) we conclude that 
  equality holds if and only if $\depth((\x_d), G_{I_1}(I_1)) \geq d-1$.
\qed

\bco
With the assumptions  as in 
Corollary~\ref{upper-bound-multiplicity} we have:
\been
\item
If  ${\displaystyle
f_{0, I_1}(I_2) 
=  e_{1}(I_{2})  -  e_{0}(I_{2}) 
- \ell \left( \f{I_2}{I_1 I_2} \right) 
- \ell \left( \f{R}{I_2} \right)
-   (d-1)~\left( \f{R}{I_1} \right)}$, then  for $0 \leq depth~G_{I_1}(I_2)$\\$ \leq d-1$,
$
\depth~F_{I_1}(I_2)=  \depth~G_{I_1}(I_2) + 1.
$

\item
If  ${\displaystyle
  f_{0, I_1}(I_2) 
  =   e_{1}(I_{2})
-      \left( 
       \f{I_1 I_2}
         {(\x_{d}) I_1 } 
       \right)
-      \sum_{n \geq 2}
       \left( 
       \f{I_1 I_2^n}
         {(\x_{d}) I_1 I_2^{n-1} } 
       \right)
+     \ell \left( 
      \f{R}{I_1} \right)}$, then  for $0 \leq depth~G_{I_1}(I_2) \leq d-2$, 
      $
\depth~F_{I_1}(I_2)=  \depth~G_{I_1}(I_2) + 1
$.
\eeen
\eco
\proof The proof follows from Corollary~\ref{upper-bound-multiplicity} and by Theorem~\ref{grade} and Theorem~\ref{depth-lemma}. 
\qed

\section{ Hilbert series for ideals of minimal and almost minimal multiplicity}
\label{hilb-fib-cone-mm}
In this section we  describe the Hilbert series of the fiber cone for ideals of minimal multiplicity and ideals of almost minimal multiplicity.

\bnot
\beqn
H_{\Fi}(i, t)
&=& \sum_{n \geq 0}H_{\Fi}(i, n)t^{n}
=  \sum_{n \geq 0}\ell \left( 
    \f{R}{I_1^i I_2^n} \right) t^n, \hspace{.5in} i=0,1.\\
H(F_{I_{1}}(I_{2}), t)
&=& \sum_{n \geq 0}H ( F_{I_{1}}(I_{2}), n)t^{n}.
\eeqn
\enot

Recall that $s(\x_d):= \min\{  n | (\x_d)I_1 I_2^n = I_1 I_2^{n+1} \}$.
\bt
\label{hilb-series-one}
Let $(R, \m)$ be a Cohen-Macaulay local ring of dimension 
$d \geq 1$.
Let $I_2$ be an ideal of minimal multiplicity with respect to $I_1$. Let $(\x_d)$ be a minimal reduction of $I_1$ and assume that $(\x_d)\subseteq I_1$.
\been
\item
${\displaystyle
    H_{\Fi}(1,t)
= \f{ \ell (R/I_1) 
- t \left[ \ell (R/I_1) -e_0(I_2) \right]  }
{(1-t)^{d}}.}$

\item
$
g_{1, I_1}(I_2) = e_0(I_2) 
- \ell \left(  \f{R}{I_1} \right).
$
For $2 \leq i \leq d$,
$g_{i, I_1}(I_2)=0.$

\item
$
{\displaystyle
H( F_{I_1}(I_2),t)
= \f{ \ell (R/I_1) - t \left[ \ell (R/I_1) -e_0(I_2) \right]  }
{(1-t)^{d+1}}
-  H_{\Fi}(0, t).
}$
\item
 $f_{0, I_1}(I_2)= e_{1}(I_{2})
    -e_{0}(I_2) 
   + \ell \left(  \f{R}{I_1} \right).$
      For $1 \leq i \leq d-1$,
 $f_{i, I_1}(I_2)=e_{i+1}(I_{2}).$
\eeen
\et
\proof 
From Theorem~\ref{homology}, 
            Lemma~\ref{mm-amm} and
and         Theorem~\ref{thm-fundamental}, 
by induction on $n$ we get
\beq
\label{alt-sum-new-g}  
     H_{\Fi}(1,n)
  =  e_0(I_2) {n+ d \choose d}
  - \left[ e_0(I_2) 
  - \ell \left( \f{R}{I_1} \right) \right]
     {n + d-1 \choose d-1}.
 \eeq
Summing over all $n \geq 0$ we get (1) and (2). 
(3) and (4) are an immediate consequence of (1), (2) and Remark~\ref{hilbert-fiber cone}. \qed

To prove Theorem~\ref{final-hilb-ser} we need the following combinatorial lemma:

\blem
\label{combin}
For all $n \geq s$ and for all $d \geq 1$ we have:
\beqn
{n-s + d-1 \choose d}
= \sum_{i=0}^{d} (-1)^i{s + 1 \choose i} {n+ d-i  \choose d-i}.
\eeqn
\elem
\proof The proof follows by induction on $d$. The case  $d=1$
can be verified easily. If $d>1$, then
\beqn
{n-s + d-1 \choose d}
&=& {n+1 -s + d-1 \choose d} - {n+1 -s + d-2 \choose d-1}\\
&=& \sum_{i=0}^{d} (-1)^i{s + 1 \choose i} {n+1+  d-i  \choose d-i}
- \sum_{i=0}^{d-1}(-1)^i{s + 1 \choose i} {n+1+  d-1-i  \choose d-1-i}\\
&=& \sum_{i=0}^{d} (-1)^i {s + 1 \choose i} {n+  d-i  \choose d-i}.
\eeqn

\bt
\label{final-hilb-ser}
Let $(R, \m)$ be a Cohen-Macaulay ring of dimension $d \geq 2$. 
Let $I_1$ and $I_2$ be an $\m$-primary ideals of $R$ and  
$(\x_{d}) \in I_2$ be a minimal reduction which is a superficial sequence for $I_1$ and $I_2$. Assume that 
$(x_d) \subseteq I_1$ and $(\x_d) \cap I_1 I_2 = I_1(\x_d)$.
\been
\item
${\displaystyle
   H_{\Fi}(1,t)
= \f{ \ell (R/I_1) 
-  t \left[ \ell (R/I_1) -e_0(I_2) \right] 
}
   {(1-t)^{d+1}}
 +  \f{t^{s+1}}{(1-t)^{d+1}}
}$

\item
$\displaystyle
{
g_{i, I_1}(I_2)=
\left\{ \begin{array}{ll}
e_0(I_2)
-   \ell \left( \f{R}{I_1} \right) 
+  s & i=1\\
{s + 1 \choose i} & i =2, \ldots, d\\
\end{array}
\right.}.
$

\item
${\displaystyle
\begin{array}{rl}
   H(F_{I_{1}}(I_{2}), t)
=& \f{ \ell (R/I_1) 
- t \left[ \ell (R/I_1) -e_0(I_2) \right] }
   {(1-t)^{d+1}}
+  \f{t^{s+1}}{(1-t)^{d+1}}
 - H_{\Fi}(0, t).
                         \end{array}
}$

\item
$\displaystyle{
   f_{i, I_1}(I_2)
= \left\{
  \begin{array}{ll}
   e_{1}(I_{2}) -e_0(I_2)
+ \ell \left( \f{R}{I_1} \right) 
-  s 
&  i=0\\
   e_{i+1}(I_{2}) -   {s+1 \choose i}    
&  i=1, \ldots,d-1
 \end{array}
\right.
}
$.
\eeen
\et
\proof Put $s = s(\x_d)$. Then 
from Theorem~\ref{homology}, 
            Lemma~\ref{mm-amm} and       Theorem~\ref{thm-fundamental}, 
\beqn
     H_{\Fi}(1,n)
  &=& \left\{ 
    \begin{array}{ll}   
          \left[ e_0(  I_2)-1 \right] {n+ d -1\choose d}
  + \ell \left( \f{R}{I_1} \right){n + d-1 \choose d-1}
  &  0 \leq n \leq s
  \\
  \hphantom{s}
 [e_0(  I_2)-1] {n+ d-1\choose d}
  + \ell \left( \f{R}{I_1} \right)
     {n + d-1 \choose d-1}
  +   {d + n - (s+1) \choose  d}
   &   n \geq s+1\\
  \end{array}
  \right.
  \eeqn
  For  all $n \geq s+1$ we have
  \beqn
   H_{\Fi}(1,n)
   &=& e_0(  I_2) {n+ d \choose d}   
   - \left[ e_0(  I_2)-1-\ell \left( \f{R}{I_1} \right) + (s+1)\right]
   + \sum_{i \geq 2} (-1)^i {s+1 \choose i} {n + d-i \choose d-i}.
   \eeqn
This proves (2). 
(3) and (4) are  an immediate consequence of (1), (2),  Remark~\ref{hilbert-fiber cone} and Lemma~\ref{combin}.  \qed

\end{document}